\documentclass[12pt]{amsart}
\usepackage{amsmath, amsfonts, amssymb, latexsym, amsthm, amscd, mathrsfs, stmaryrd} 
\usepackage[all]{xy}
\usepackage[makeroom]{cancel}
\usepackage{yhmath}
\usepackage{hyperref}

\theoremstyle{definition}
\newtheorem{rem}[subsubsection]{Remark}
\theoremstyle{plain}
\newtheorem{prop}[subsubsection]{Proposition}
\newtheorem{thm}[subsubsection]{Theorem}
\newtheorem{lem}[subsubsection]{Lemma}

\newcommand{\mbf}{\mathbf}
\newcommand{\mbb}{\mathbb}
\newcommand{\mrm}{\mathrm}

\newcommand{\E}{\breve E}
\newcommand{\F}{\breve F}

\newcommand{\K}{\breve K}

\renewcommand{\S}{\mbf S}

\newcommand{\U}{\mbf U}
\newcommand{\V}{\mbf V}

\renewcommand{\k}{\mbf k}

\newcommand{\hsl}{\widehat{\mathfrak{sl}}}
\newcommand{\hgl}{\widehat{\mathfrak{gl}}}
\newcommand{\ve}{\varepsilon}

\usepackage{color}
\newcommand{\nc}{\newcommand}
\nc{\redtext}[1]{\textcolor{red}{#1}}
\nc{\bluetext}[1]{\textcolor{blue}{#1}}
\nc{\greentext}[1]{\textcolor{green}{#1}}
\nc{\yl}[1]{\redtext{ #1}}
\nc{\zb}[1]{\redtext{From zb: #1}}

\title[Embeddings among quantum affine $\mathfrak{sl}_n$]{Embeddings among  quantum affine $\mathfrak{sl}_n$}

\author{Yiqiang Li$^{\dagger}$} \thanks{$\dagger$Partially supported by NSF DMS 1801915}
\address{Department of Mathematics\\ University at Buffalo\\ the State University of New York  \\Buffalo, NY 14260}
\email{yiqiang@buffalo.edu}

\date{\today }
\keywords{Quantum affine $\mathfrak{sl}_n$, embeddings}
\subjclass{17B37}

\begin{document}

\begin{abstract}
We establish an explicit embedding of a quantum affine $\mathfrak{sl}_n$ into a quantum affine $\mathfrak{sl}_{n+1}$. This embedding serves as a common generalization of 
two natural, but seemingly unrelated,  embeddings, one on the quantum affine Schur algebra level and  the other on the non-quantum level. 
The embedding on the quantum affine Schur algebras is used extensively in the analysis of canonical bases of quantum affine $\mathfrak{sl}_n$
and $\mathfrak{gl}_n$.
The embedding on the non-quantum level is used crucially in a work of Riche and Williamson on the study of modular representation theory
of general linear groups over a finite field. The same embedding is also used in a work of Maksimau on the categorical representations
of affine general linear algebras. 
We  further provide a more natural compatibility statement of the embedding on the idempotent version with that on the quantum affine Schur algebra level.
A $\hgl_n$-variant of the embedding is also established. 
\end{abstract}

\maketitle

\section*{Introduction}
Consider the following rule
\begin{align}
\label{rule}
\begin{bmatrix}
a & b\\
c & d
\end{bmatrix}
\mapsto 
\begin{bmatrix}
a& 0 & b\\
0 & 0 & 0 \\
c & 0 & d
\end{bmatrix}.
\end{align}
It defines an embedding of $\mathfrak{sl}_2(R)$  into $\mathfrak{sl}_3(R)$  over any commutative ring $R$ where $a, b, c$ and $d$ belong to. 
Further, by regarding $a, b, c$ and $d$ as block matrices,
the embedding generalizes to  a natural embedding of $\mathfrak{sl}_n(R)$ into $\mathfrak{sl}_{n+1}(R)$.
When $R$ is the local field $\mbb C((t))$,
such an embedding
plays a key role in 
the study of categorical representations of affine general Lie algebras in Maksimau's work ~\cite{M15, M18}.
An affine $\mathfrak{gl}$ variant of the embedding is further used  in the study of modular representations of general linear groups over finite fields in a recent work ~\cite{RW18} of Riche and Williamson.

Now consider the $n$-step affine flag variety  $\mathfrak F_{n, d}$ associated with $\mathrm{GL}_d(\mbb C)$. 
It is well known that quantum affine Schur algebras  $\S_{n, d}$ admit a geometric realization as the convolution algebra on   $\mathfrak F_{n, d}\times \mathfrak F_{n, d}$. 
The ind-variety $\mathfrak F_{n, d}$ parametrizes lattice chains in a $d$-dimensional vector space over $\mbb C((t))$. 
The operation of adding an extra copy of the lattice in a prescribed step in lattice chains defines a natural embedding $\S_{n, d}\subset \S_{n+1, d}$. 
Such a natural embedding of quantum affine Schur algebras is used crucially in the study of canonical bases and multiplication formulas of quantum affine 
$ \mathfrak{sl}_n(\mbb C)$ and $\mathfrak{gl}_n(\mbb C)$ in~\cite{FLLLW, LS}. 
We refer to~\cite{DDF} for a Hecke-algebra approach to, and  further applications of, this embedding.

So we have a natural embedding of affine $\mathfrak{sl}_n$ into affine $\mathfrak{sl}_{n+1}$ on the one hand and a natural embedding $\S_{n,d}\subset \S_{n+1,d}$ on the other hand. 
Note that $\S_{n,d}$ receives an algebra homomorphism from quantum affine $\mathfrak{sl}_n$. To this end, 
it is desirable to understand the relationship between the two kinds of seemingly unrelated embeddings. 
In this paper, we establish an explicit embedding of quantum affine $\mathfrak{sl}_n$ into its higher rank. 
It is a natural quantization of the embedding on affine $\mathfrak{sl}_n$ defined by (\ref{rule}) and used in the works~\cite{RW18, M15, M18}. 
We further provide a compatibility of the embedding with that on the quantum affine Schur algebra level. It is somewhat unnatural. However when consider the embedding on Lusztig's modified quantum affine $\mathfrak{sl}_n$,  it becomes natural again.  So the embeddings established in this paper can be regarded as a common generalization of 
the previous two kinds of embeddings. 
It is worthwhile to point out that the quantization  is not unique, depending on a parameter $\ve\in \{ \pm 1\}$. 
In an upcoming paper~\cite{Li}, the author will extend the embeddings established in this paper to a much broader setting. 

Note that if one replaces the zero in the center of the right matrix by
$- \mrm{trace} (a) - \mrm{trace}(d)$
 in the assignment (\ref{rule}), it induces an embedding $\mathfrak{gl}_n(R)\subseteq \mathfrak{sl}_{n+1}(R)$. 
As a quantization of this embedding, we further provide an embedding of quantum affine $\mathfrak{gl}_n$ in the sense of Green~\cite{G99}  into quantum affine $\mathfrak{sl}_{n+1}$ of level $1$.
There is a second version of quantum affine $\mathfrak{gl}_n$ studied by Du and Fu~\cite{DF15}. 
It is very interesting to see if there is  a similar embedding  as a common generalization of the embeddings on the non-quantum level and on the quantum affine Schur algebra level.

 
 Let $\theta_n$ be an involution of $\mathfrak{sl}_n(R)$ by sending a matrix to the matrix whose $(i, j)$ entry is the $(n+1-i, n+1-j)$ entry of the original matrix. 
The pair $(\mathfrak{sl}_n(R), \mathfrak{sl}_{n}(R)^{\theta_n})$ is a quasi-split symmetric pair.
The embedding defined in (\ref{rule}) with appropriate block sizes of $a$ and $d$ is compatible with the involutions $\theta_{n}$ and $\theta_{n+1}$ on $\mathfrak{sl}_{n}(R)$ and $\mathfrak{sl}_{n+1}(R)$ respectively. Thus it induces an
embedding on the fixed-point set: $\mathfrak{sl}_n(R)^{\theta_{n}}\subset \mathfrak{sl}_{n+1}(R)^{\theta_{n+1}}$.  
On the geometry side, there is a quantum affine   Schur algebras $\S_{n, d}^{\imath}$ defined as the convolution algebras of $n$-step affine isotropic flag varieties 
 in the works~\cite{FLLLW} (see also~\cite{BKLW}). 
 These algebras are homomorphic images of a quantum version of $\mathfrak{sl}^{\theta_n}_n$, i.e., the coideal subalgebras in an affine quantum symmetric pair of quasi-split   type $A$. There are  natural embeddings $\S_{n,d}^{\imath} \subseteq \S_{n+1, d}^{\imath}$, which play a key role in {\it loc. cit}.
It is a natural question to see if there exists an embedding for these types of  coideal subalgebras as a common generalization of the embedding 
$\mathfrak{sl}_n(R)^{\theta_{n}}\subset \mathfrak{sl}_{n+1}(R)^{\theta_{n+1}}$ and $\mathfrak{sl}_n(R)^{\theta_{n}}\subset \mathfrak{sl}_{n+1}(R)^{\theta_{n+1}}$.

 The embedding established in this paper certainly calls for a further investigation of a  more intrinsic connection between the above lines of research.

\tableofcontents

\section{Preliminaries}

This section recalls the basic definitions of quantum affine $\mathfrak{sl}_n$ and its Schur quotients. 

\subsection{Definition}
Let $\mbb Q(v)$ be the field of rational functions with the variable $v$.
For each integer $a\in \mbb N$, we define
\[
[a]= \frac{v^a-v^{-a}}{v-v^{-1}}, [a]^! = [a] [a-1] \cdots [1].
\]
If $x$ is an element in an associative algebra over $\mbb Q(v)$, we write
\[
x^{(a)} = x^a/[a]^!
\]

Let $I_n=\mbb Z/n\mbb Z$.
 If there is no danger of confusion, we write elements in $I_n$ by $ 0, 1, \cdots, n-1$. 
Recall that the Cartan matrix  of affine type $A^{(1)}_{n-1}$ is the $n$ by $n$ matrix 
$C=(c_{ij})_{i, j\in I_n} $ such that 
$c_{ij} = 2 \delta_{i,j} -\delta_{i, j+1} - \delta_{i, j-1}$.
Let $\U(\hsl_n)$ be the quantum affine $\hsl_n$ associated with $C$. It is an associative algebra over $\mbb Q(v)$
and it admits a generator-and-relation presentation. 
Precisely, the generators are  
standard Chevalley generators $ E_i, F_i,$ and $ K_i^{\pm 1}$ for all  $ i\in I_n$ and the  defining relations are given as follows. 
\begin{align}
\label{R1}
&K_{i}   K_{j}  =    K_{j}   K_{i},\
 K_{i}  K^{-1}_{i}  =1=K^{-1}_i K_i,  \tag{R1}\\
 \label{R2}
  &K_{i}  E_{j} 
 = v^{c_{ij}}  E_{j}  K_i, \
  K_{i}  F_{j}   
= v^{-c_{ij} }  F_{j}  K_i, \tag{R2}\\
 \label{R3}
 &E_{i}  F_{j} -  F_{j}  E_{i} = \delta_{i,j} \frac{  K_i
 -  K^{-1}_{i} }{v-v^{-1}},        \tag{R3}  \\
 \label{R4}
 & \sum_{p=0}^{1-c_{ij}} (-1)^p E_i^{(p)} E_j E_i^{(1-c_{ij} -p)}=0,\tag{R4} \\
 \label{R5}
 &\sum_{p=0}^{1-c_{ij}} (-1)^p F_i^{(p)} F_j F_i^{(1-c_{ij} -p)}=0, \forall i, j \in I_n . \tag{R5}
 \end{align}

The algebra $\U(\hsl_n)$ can be endowed with a Hopf algebra structure with the comultiplication $\Delta$ defined by
\begin{align*}
\Delta(E_i) = E_i\otimes 1 + K_i \otimes E_i,
\Delta(F_i) = F_i\otimes K^{-1}_i + 1\otimes F_i,
\Delta(K_i) = K_i \otimes K_i,\quad \forall i \in I_n.
\end{align*}

\subsection{Idempotent form}

\label{idem}

Let $n$ be a positive integer greater than $1$. 
Let $\mathfrak S_n$ be the set of all sequences $(a_i)_{i\in \mbb Z}$ such that $a_i\in \mbb Z$ and $a_i=a_{i+n}$ for all $i$.
Note that any sequence $(a_i)_{i\in \mbb Z}$ in $\mathfrak S_n$ is completely determined by the entries
$(a_1,\cdots, a_n)$, which we shall use to denote the sequence itself.
The set $\mathfrak S_n$   carries a natural abelian group structure with the component-wise addition.
It is naturally isomorphic to $\mbb Z^n$.
Let $Y_n $ be the subset of $\mathfrak S_n$ consisting of all sequences $(a_i)_{i\in \mbb Z}$ such that $\sum_{1\leq i\leq n} a_i=0$.
Let $\mbf 1$ be the sequence in $\mathfrak S_n$ such that each entry is $1$. Let 
$X_n= \mathfrak S_n/ \langle \mbf 1\rangle $ be the quotient of $\mathfrak S_n$ by the subgroup generated by $\mbf 1$.
Let $\langle -, -\rangle: Y_n\times X_n\to \mbb Z$ be the perfect pairing defined by
$\langle \mbf b, \bar{\mbf a} \rangle=\sum_{1\leq i\leq n} b_ia_i$ for all $\mbf b\in Y_n$ and $\bar{\mbf a} \in X_n$.
Let $I_n=\mbb Z/n\mbb Z$. 
There is an embedding $I_n\to Y_n$ defined by $i\mapsto \alpha_i^{\vee}$ where the $i$ (resp. $i+1$) entry of $\alpha_i^{\vee}$ is 
$1$ (resp. $-1$) mod $n$ and $0$ otherwise.
Similarly there is an embedding $I_n\to X_n$ defined by $i\mapsto \alpha_i$ where $\alpha_i$ is the coset of $\alpha^{\vee}_i$ in $X_n$. 
Let $i\cdot j=\langle \alpha^{\vee}_i, \alpha_j\rangle$ for all $i,j\in I_n$. Then $(I_n, \cdot)$ is the Cartan datum of affine type $A_n$. 
Further,  $(Y_n, X_n, \langle-,-\rangle)$ is a root datum of $(I_n, \cdot)$.
The root datum is neither $X$-regular nor $Y$-regular. 

Let $\dot \U_n$ be the modified quantum group associated with the root datum $(Y_n, X_n)$ (\cite{L10}). 
The algebra $\dot \U_n$ is an associative algebra over $\mbb Q(v)$ without unit. 
Instead it has a collection of idempotents $1_{\lambda}$ for $\lambda\in X_n$. 
It carries a natural $\U(\hsl_n)$-bimodule structure generated by  the $1_{\lambda}$'s.

\subsection{Quantum affine Schur algebra}
\label{tensor}

Let $\V_n$ be the vector space over $\mbb Q(v)$ spanned by the symbols $u_k$ for all $k\in \mbb Z$.
To avoid confusion, we write $\bar k$ for  the congruence class of $k$ in $I_n=\mbb Z/n\mbb Z$.
We define a $\U(\hsl_n)$-module structure on $\V_n$ by the following rules: for all $\bar i\in I_n$, $k\in \mbb Z$, 
\begin{align*}
E_{\bar i}. u_k  = \delta_{\bar k,\overline{ i+1}} u_{k-1},\
F_{\bar i} . u_k  =\delta_{\bar k,  \bar i} u_{k+1},\
K_{\bar i}. u_k  = v^{\delta_{\bar k, \bar  i} - \delta_{\bar k, \overline{i+1}}} u_k.
\end{align*}
Via the comultiplication $\Delta$, there is a $\U(\hsl_n)$-module structure on the tensor product 
$\V_{n}^{\otimes d} =\V_n\otimes \cdots \otimes \V_n$ of $d$ copies of $\V_n$. 
There exists a natural action of the affine Hecke algebra $H_d$ of type $A$ on $\V_n^{\otimes d}$.
The quantum affine Schur algebra is defined to be the centralizer algebra
$\S_{n, d} = \mrm{End}_{H_d}(\V^{\otimes d})$. 
We shall not recall the precise definition of the $H_d$-action as it plays no role in later analyses, instead we shall recall 
a geometric definition of $\S_{n, d}$ in the following. 
Since the two actions on $\V^{\otimes d}$ are commutative, there exists an algebra homomorphism 
$\U(\hsl_n)\to \S_{n, d}$. The morphism is not surjective in general. Let 
$\S'_{n,d}$ be its image,  a subalgebra in  
$\S_{n,d}$. 
The algebra $\S'_{n,d}$ was first studied by Lusztig. 
In particular, we have a surjective algebra homomorphism
\begin{align}
\label{pi}
\pi_{n,d}: \U(\hsl_n)\to \S'_{n,d}.
\end{align} 
It is well known from~\cite{L00} that there exists an algebra homomorphism $\phi_{d+n, d}: \S'_{n,d+n} \to \S'_{n,d}$ such that 
$\pi_{n, d} = \phi_{d+n,d} \pi_{d+n,d}$. Moreover $\U(\hsl_n)$ is in the inverse limit of the inverse system
$\{ \S_{n, d} , \phi_{d+n,d}\}_{d\in \mbb N}$. 

Let $\k$ be a finite field of $q$ elements. 
Let $\k((t))$ be the field of formal Laurent polynomials. 
Let $\k[[t]]$ be the subring of $\k((t))$ consisting of all formal power series. 
Consider a $\k((t))$-vector space $V$  of dimension $d$. 
A  free $\k[[t]]$-module $L$ of rank $d$ such that $\k((t)) \otimes_{\k[[t]]} L= V$ is called a lattice of $V$. 
The collection of chains of lattices $L_\bullet= (L_1\subseteq L_2\subseteq \cdots \subseteq L_n\subset t^{-1} L_1)$ will be denoted by
$\mathfrak F_{n, d}$, the $n$-step affine flag variety. 
When specialized at $v= q^{1/2}$, the algebra $\S_{n, d}$ is isomorphic to the convolution algebra on
$\mathfrak F_{n, d}\times \mathfrak F_{n, d}$.

Let $\Lambda_{n,d}$ be the subset in $\mathfrak S_n$ of all sequence $(a_i)$ such that $\sum_{1\leq i\leq n} a_i=d$.
For each $(a_i) \in \Lambda_{n,d}$, define $\mathfrak F_{(a_i)}$ to be the subset of $\mathfrak F_{n, d}$ 
consisting of all lattice chains $L_\bullet$ such that $\dim_{\k}  L_i/L_{i-1}= a_i$ for all $1\leq i\leq n$.
For each $(a_i)\in \Lambda_{n, d}$, there is an idempotent $1_{(a_i)}$ in $\S_{n, d}$ such that, when specialized to
$v=q^{1/2}$, it is the characteristic function of the diagonal of $\mathfrak F_{(a_i)}\times \mathfrak F_{(a_i)}$.

\section{The embedding $\Phi_{r,\ve}$}

In this section, we present an explicit embedding of the quantum affine $\mathfrak{sl}_n$ to its higher rank 
and a compatibility with the natural embedding on the Schur algebra level.

\subsection{The statement}
Let us fix forever an integer $r\in [0, n-1]$ and $\ve\in \{\pm 1\}$. Consider  $\U(\hsl_{n+1})$. To avoid confusion, we denote the  Chevalley generators
in $\U(\hsl_{n+1})$ by
$ \E_i, \F_i, \K_i^{\pm 1}$ for $i\in I_{n+1}$.
Consider the following elements in $\U(\hsl_{n+1})$.
For all $i\in I_n$, we define

\begin{align}
\label{efk}
\begin{cases}
\begin{split}
e_{i} & =
\begin{cases}
\E_i & \mbox{if} \ i= 0, \cdots, r-1,\\
\E_i \E_{i+1} - v^{\ve} \E_{i+1} \E_i & \mbox{if}\ i= r,\\
\E_{i+1} & \mbox{if} \ i= r+1, \cdots, n-1.
\end{cases}\\
f_{i} & =
\begin{cases}
\F_i &  \mbox{if} \ i= 0,\cdots, r-1,\\
\F_{i+1} \F_i - v^{-\ve} \F_i\F_{i+1} & \mbox{if} \ i=r,\\
\F_{i+1} & \mbox{if} \ i =r+1,\cdots, n-1. 
\end{cases}\\
k_i & =
\begin{cases}
\K_i &\mbox{if} \ i=0,\cdots, r-1,\\
\K_i\K_{i+1} & \mbox{if} \ i= r,\\
\K_{i+1}  & \mbox{if} \ i=r+1, \cdots, n-1.
\end{cases}
\end{split}
\end{cases}
\end{align}

Recall that for any $x\in \U(\hsl_{n+1})$ we write $x^{(a)}=x^a/[a]^!$. 
For any $p\in \mbb N$, we have
\begin{align}
\label{e-p}
e_r^{(p)} = \sum_{j=0}^p (-1)^j v^{\ve j} \E_{r+1}^{(j)} \E^{(p)}_r \E_{r+1}^{(p-j)}, \\
f_r^{(p)} =  \sum_{j=0}^p (-1)^j v^{-\ve j} \F_{r}^{(j)} \F^{(p)}_{r+1} \F_{r}^{(p-j)}. 
\end{align}
The equality (\ref{e-p}) is ~\cite[Lemma 42.1.2 (c)]{L10}.
We have the following.

\begin{thm}
\label{Phi}
There is  an injective homomorphism of associative algebras
\begin{align}
\begin{split}
\Phi_{r,\ve} :  \U(\hsl_n) \to \U(\hsl_{n+1})  \ \mbox{defined by}\
 E_i\mapsto e_i, F_i\mapsto f_i, K^{\pm 1}_i\mapsto k_i^{\pm 1},\ \forall  i\in I_n.
\end{split}
\end{align}
\end{thm}

\begin{proof}

We shall show the existence of  $\Phi_{r,\ve}$ in Section~\ref{existence}: 
for $n>2$, it is given in Section~\ref{Proof-1} and  for $n=2$ it is given in Section~\ref{Proof-2}.

Assume now that the  algebra homomorphism $\Phi_{r,\ve}$ is well-defined.  
We shall show that the map $\Phi_{r,\ve}$ is injective.
Let $\U_{r,\ve}$ be the subalgebra of $\U(\hsl_{n+1})$ generated by 
$e_i, f_i, k^{\pm 1}_i$ for all $i\in I_n$. 
We consider the module $\V_{n+1}$ of $\U(\hsl_{n+1})$.
It restricts   automatically  to a $\U_{r,\ve}$-module. 
To avoid confusion, we write its basis element by $\breve u_k$ for all $k\in \mbb Z$.
We define an embedding $\V_n\to \V_{n+1}$ by
$$
u_k \mapsto 
\begin{cases}
\breve u_{k+ \lfloor k/n \rfloor} &\mbox{if} \  \bar k \in \{ \overline {0}, \overline {1}, \cdots, \overline{ r-1} \} \subseteq I_n,\\
\breve u_{k+\lceil k/n\rceil} & \mbox{if} \ \bar k \in \{ \overline{ r},\cdots, \overline{n-1}\} \subseteq I_n.
\end{cases}
$$
Let $\mbf W_n$ be the image of the embedding $\V_n\to \V_{n+1}$. 
One can check that $\mbf W_n$ is a $\U_{r,\ve}$-submodule of $\V_n$, and moreover
the actions  of the generators in  $\U_{r,\ve}$ on $\mbf W_n$ is compatible with the action of
the generators in $\U(\hsl_n)$ under the isomorphism $\V_n\to \mbf W_n$.
Observe that 
\begin{align}
\label{D-e}
\Delta (e_r) & = e_r\otimes 1 + k_r \otimes e_r + (v^{-1} - v^\ve) \E_{r+1} \K_r \otimes \E_r +(v-v^{\ve}) \K_{r+1} \E_r \otimes \E_{r+1}, \\
\Delta(f_r) & =  f_r\otimes k^{-1}_r + 1\otimes f_r + (v-v^{-\ve}) \F_r \otimes \K^{-1}_r \F_{r+1} + (v^{-1}-v^{-\ve}) \F_{r+1} \otimes \F_r \K^{-1}_{r+1}.
\end{align}
This indicates that when restricting to $\mbf W_n^{\otimes d}$, as a summand of $\V_{n+1}^{\otimes d}$, the last two operators
in $\Delta(e_r)$ and $\Delta(f_r)$ are zero. 
Thus $\mbf W_n^{\otimes d}$ is a $\U_{r,\ve}$-module and thanks to (\ref{D-e})
the actions of $e_r$ and $f_r$
on $\mbf W_n^{\otimes d}$ are compatible with the actions of $E_r, F_r$ on $\V_n^{\otimes d}$ via the isomorphism $\V_n\to \mbf W_n$.
Identifying $\V_n $ with $\mbf W_n$, we see that there exists an algebra homomorphism $\pi'_{n,d}: \U_{r,\ve} \to \S'_{n,d}$ such that 
$\pi'_{n,d} (e_i)= \pi_{n,d} (E_i)$, $\pi'_{n,d}(f_i) = \pi_{n,d} (F_i) $ and $\pi'_{n,d} (k_i)=\pi_{n,d}(K_i)$ for all $i\in I_n$.
Since $\U(\hsl_n)$ is in the inverse limit of an inverse system $\{ \S'_{n,d}, \S'_{n,d+n}\to \S'_{n,d}\}_{d\in \mbb N}$, by the universality there exists a unique algebra homomorphism $\Psi_{r,\ve}: \U_{r,\ve} \to \U(\hsl_n)$ such that
$\Psi_{r,\ve} (e_i)=E_i$, $\Psi_{r,\ve}(f_i)= F_i$ and $\Psi_{r,\ve}(k_i^{\pm 1})= K^{\pm 1}_i$ for all $i\in I_n$. This is the inverse of the map
$\Phi'_{r,\ve}: \U(\hsl_n) \to \U_{r,\ve}$ induced by $\Phi_{r,\ve}$. So $\Phi'_{r,\ve}$ is an isomorphism, and in other words, $\Phi_{r,\ve}$ is injective.
The theorem is thus proved. 
\end{proof}

We end this section with a remark.

\begin{rem}
\begin{enumerate}
\item The assumption that $r\in [0, n-1]$ is not essential. One can define the morphism $\Phi_{r,\ve}$ when $r=n$ as well in a similar manner.

\item The $v=1$ version of the embedding $\V_n\to \V_{n+1}$ and (\ref{efk})  first appeared in~\cite{M18}. 

\item The analysis in the Proof of Theorem~\ref{Phi} yields the following commutative diagram of linear maps. 
\[
\begin{CD}
\U(\hsl_n) @> \Phi_{r,\ve} >> \U_{r,\ve} \\
@V\pi_{n,d} VV @VV\tilde \pi_{n+1,d} V\\
\S'_{n,d} @< \pi << \S'_{n+1,d}
\end{CD}
\]
where 
$\tilde \pi_{n+1,d}$ is the restriction of the map $\pi_{n+1,d}$ in (\ref{pi}) to $\U_{r,\ve}$ and 
$\pi $ is the projection of a linear endomorphism of $\V_{n+1}^{\otimes d}$ to a linear endomorphism 
of $\V_n^{\otimes d}$ under the embedding $\V_n\to \V_{n+1}$. 
(Note that $\pi$ is not an algebra homomorphism.)
A more natural compatibility will be given in the following section. 
The injectivity of $\Phi_{r,\ve}$ can be proved by exploring the above commutative diagram as follows.
Suppose that $x\in \U(\hsl_n)$ is in the kernel of $\Phi_{r,\ve}$. 
Then by the above commutative diagram, we have $\pi_{n,d} (x)= \pi \pi'_{n+1,d} \Phi_{r,\ve} (x) =0$ for all $d$.
It is well-known that if $\pi_{n,d}(x)=0$ for all $d$, then  $x=0$. Therefore $\Phi_{r,\ve}$ is injective. 

\end{enumerate}
\end{rem}

\subsection{A variant}

In~\cite{G99}, Green defines a version of quantum affine $\mathfrak{gl}_n$. This is a variant of $\U(\hsl_n)$.
 More precisely, it is defined as an associative algebra $\U(\hgl_n)$ over $\mbb Q(v)$, which has
generators $E_i, F_i, L^{\pm 1}_i$ for $i\in I_n$ and which subject to the following defining relations.
\begin{align}
\label{gl-def}
\begin{cases}
&L_i L_j=L_j L_i , L_i L_i^{-1}=1,\\
&L_i E_j = v^{\delta_{i,j}  -\delta_{i-1, j}} E_j L_i, \\
&L_i F_j = v^{-\delta_{i,j}  + \delta_{i-1, j}} F_j L_i,\\
&E_i F_j -F_j E_i =\delta_{i,j} \frac{L_i L_{i+1}^{-1} - L^{-1}_i L_{i+1}}{v-v^{-1}},\\
& \sum_{p=0}^{1-c_{ij}} (-1)^p E_i^{(p)} E_j E_i^{(1-c_{ij} -p)}=0, \\
&\sum_{p=0}^{1-c_{ij}} (-1)^p F_i^{(p)} F_j F_i^{(1-c_{ij} -p)}=0, \forall i, j \in I_n . 
\end{cases}
\end{align}

Note that in~\cite{G99}, there is an assumption that $n\geq 3$. We do not need this assumption. 

Recall that, to avoid confusion, a superscript $\breve \empty$ is put on the Chevalley generators of $\U(\hsl_{n+1})$.
Recall that we fix an integer $r\in [0, n-1]$. Consider the following element in $\U(\hsl_{n+1})$.
\begin{align}
l_i =
\begin{cases}
\breve K_i \breve K_{i+1} \cdots\breve K_r & 1\leq i\leq r,\\
\breve K^{-1}_{r+1} \cdots \breve K^{-1}_{i} & r+1\leq i\leq n. 
\end{cases}
\end{align}

Let $\U(\hsl_{n+1})_1\equiv \U(\hsl_{n+1}) / \langle \prod_{i=1}^{n+1} \breve K_i - 1\rangle $ be the quotient algebra of $\U(\hsl_{n+1})$ by the two-sided ideal 
generated by $\prod_{i=1}^{n+1} \breve K_i - 1$. 
This is the quantum affine $\hsl_{n+1}$ at level $1$. 
We still use the same notations for the images of the Chevalley generators of $\U(\hsl_{n+1})$ under the canonical projection map. 
We have the following $\widehat{\mathfrak{gl}}_n$-variant of Theorem~\ref{Phi}

\begin{prop}
\label{variant}
There is an embedding $\Phi_{r, \ve; 1}: \U(\hgl_n) \to \U(\hsl_{n+1})_1$ defined by 
$E_i\mapsto e_i, F_i \mapsto f_i, L_i\mapsto l_i$ for all $i\in I_n$, where $e_i, f_i$ are in (\ref{efk}). 
\end{prop}

The proof is given in Section~\ref{Proof-variant}.

\subsection{Compatibility with quantum affine Schur algebras}
\label{Schur}

Recall the root datum $(Y_n, X_n)$ from Section~\ref{idem}. 
Define a morphism    $$\phi=(f,g): (Y_n,X_n)\to (Y_{n+1},X_{n+1})$$ where
$f: Y_n\to Y_{n+1}$ sends a sequence $(a_1,\cdots, a_n)$ to $(a_1, \cdots, a_{r}, 0, a_{r+1} , \cdots, a_n)$, and
$g: X_{n+1}\to X_n$ sends a coset of $(a_1,\cdots, a_{n+1})$ to the coset of $(a_1- a_{r+1}, \cdots, a_{i_r}-a_{r+1} , a_{r+2}-a_{r+1},\cdots, a_{n+1} - a_{r+1})$. 
The root datum $(Y_{n+1}, X_{n+1})$ can be regarded as a root datum of affine type $A_n$, instead of $A_{n+1}$, 
in an appropriate way such that $\phi$ is a morphism of root data of affine type $A_n$.
Note that $f$ extends naturally to a map $f: \mathfrak S_n\to \mathfrak S_{n+1}$.


There is an embedding $\Lambda_{n, d} \to X_n$  defined by $\lambda \mapsto \bar \lambda$
as the composition of the embedding $\Lambda_{n, d} \to \mathfrak S_n$ and the quotient map $\mathfrak S_n \to X_n$. 
For $\bar \lambda \in X_n - \Lambda_{n, d}$, then there is a $d' \neq d$ such that $\lambda \in \Lambda_{n, d'}$. 
Clearly different representatives of $\bar \lambda $ yield different $d'$. 
We fix once and for all a representative for each class in $X_n$ so that if $\bar \lambda \cap \Lambda_{n, d}\neq \O$ then $\lambda \in \Lambda_{n, d}$.
Under this assumption, we can define a map $f_d: X_n \to X_{n+1}$ by sending $\bar \lambda \to \overline{f(\lambda)}$. 
Clearly $g f_d (\bar \lambda ) =\bar \lambda$.
In light of Theorem~\ref{Phi}, we thus have an algebra embedding 
\begin{align}
\label{phi}
\phi_{d,\ve}: \dot \U_n \to \dot \U_{n+1}
\end{align}
defined  by $1_{\bar \lambda} \mapsto 1_{\overline{f(\lambda)}}$, $E_i 1_{\bar \lambda} \mapsto e_i 1_{\overline{f(\lambda)}}$ and 
$F_i 1_{\bar \lambda} \mapsto f_i 1_{\overline{f(\lambda)}}$ for all $i\in I_n$ and $\lambda \in \mathfrak S_n$.

Recall from~\ref{tensor} the quantum affine Schur algebra $\S_{n,d}$. It has  idempotents $1_{\lambda}$   for $\lambda \in \Lambda_{n, d}$. 
It is known that there is an algebra homomorphism $\dot \pi_{n, d}: \dot \U_n\to \S_{n,d}$ such that 
$\dot \pi_{n, d} (1_{\bar \lambda}) = 1_{\lambda}$ if $ \lambda \in \Lambda_{n, d}$, and
$\dot \pi_{n, d} (1_{\bar \lambda})=0$ otherwise.
One the other hand, the subset $\mathfrak F_{n+1, d}|_r$ of $\mathfrak F_{n+1, d}$ consisting of all lattice chains $L_\bullet$ such that $L_r=L_{r+1}$ is naturally in bijection with $\mathfrak F_{n,d}$. Via this bijection, 
there is an algebra homomorphism $\sigma_d: \S_{n,d} \to \S_{n+1,d}$ such that $\sigma_d ( 1_{\lambda}) = 1_{f(\lambda)}$ 
for all $\lambda\in \Lambda_{n, d}$ (see e.g., ~\cite{LS}). 
We are ready to state the compatibility of the embeddings on $\dot \U_n$ and $\S_{n, d}$.

\begin{prop}
\label{Psi-Schur}
Let $\ve\in \{ \pm 1\}$. 
The following  diagram is commutative.
\begin{align}
\label{Psi-Schur-a}
\begin{split}
\xymatrix{
\dot \U_n \ar[rr]^{\dot \phi_{d,\ve}} \ar[d]_{\pi_{n, d}}  &&  \dot \U_{n+1} \ar[d]^{\dot \pi_{n+1,d}}  \\
\S_{n, d} \ar[rr]_{\sigma_d} && \S_{n+1, d}
}
\end{split}
\end{align}
where $\phi_{d,\ve}$ is from (\ref{phi}).
\end{prop}

\begin{proof}
By definition, the two compositions $\sigma_d\circ \dot \pi_{n, d} $ and $\dot \pi_{n+1,d} \circ  \phi_{d,\ve}$ coincide when evaluating at the idempotents
$1_{\lambda}$ for all $\lambda \in X_n$. Now the maps are all $\U_n$-module homomorphisms and hence the two compositions coincide when
evaluating at any element in $\dot \U_n$. Hence the diagram must be commutative. The proposition follows.
\end{proof}

\section{Existence of $\Phi_{r,\ve}$}

\label{existence}

This section is devoted to the proof of the existence of $\Phi_{r,\ve}$. 
In Section~\ref{Proof-1}, we present the proof when $n>2$ and in Section~\ref{Proof-2}, we present the proof for the $n=2$ case. 

\subsection{Existence of $\Phi_{r,\ve}$ in the $n>2$ case}
\label{Proof-1}

In this section we assume that $n>2$. We shall show that $\Phi_{r,\ve}$ is an algebra homomorphism, i.e., the set
$\{e_i, f_i, k^{\pm 1}_i| i\in I_n\}$ satisfies the defining relations of $\U(\hsl_n)$.

It is easy to see that the elements $k^{\pm 1}_i$ satisfy the relation (\ref{R1}).

Next, we show that the pair $(k_i, e_j)$ satisfies the defining condition in (\ref{R2}). All are a consequence of the definition, except the case 
$(k_i, e_i)$, $(k_{i\pm 1}, e_i)$ for $i =r$. Observe that 
\[
\K_{i-1} e_i = v^{-1} e_i \K_{i-1} , \K_i e_i = v e_i \K_i, \K_{i+ 1} e_i = v e_i \K_{i+ 1},\K_{i+2} e_i = v^{-1} e_i \K_{i+2}  \ \mbox{if} \ i =r. 
\]
So we have 
\begin{align*}
k_i e_i = \K_i \K_{i+1} e_i = v^2 e_i \K_i\K_{i+1} = v^2 e_i k_i.\
k_{i\pm 1}  e_i = v^{-1} e_i k_{i\pm 1}.
\end{align*}
Thus  the pair $(k_i, e_j)$ satisfies the relation in (\ref{R2}).
Similarly the pair $(k_i, f_j)$ satisfies the relation in (\ref{R2}) because we have 
\[
\K_{i-1} f_i = v f_i \K_{i-1} , \K_i f_i = v^{-1} f_i \K_i, \K_{i+ 1} f_i = v^{-1} f_i \K_{i+ 1},\K_{i+2} f_i = v f_i \K_{i+2}  \ \mbox{if} \ i =r. 
\]

Next, we show that the triple $(e_i, f_j, k^{\pm 1}_i)$ satisfies the commutator relation in (\ref{R3}).
The relation is satisfied automatically if $i, j \neq r$. 
Assume that $i=r, j\neq r$, and $j< r$.
Then we have
\begin{align*}
e_i f_j-f_je_i 
&= (\E_i\E_{i+1} - v^{\ve} \E_{i+1} \E_{i} ) \F_j - \F_j ( \E_i \E_{i+1} - v^{\ve} \E_{i+1} \E_i) \\
& = \E_i \E_{i+1} \F_j - \F_j \E_i \E_{i+1} - v^{\ve} ( \E_{i+1} \E_i \F_j - \F_j \E_{i+1} \E_i)\\
&=(\E_i\F_j  - \F_j \E_i) \E_{i+1} - v^{\ve}  \E_{i+1} (\E_i \F_j -\F_j  \E_i)=
0. 
\end{align*} 
So the commutator relation holds for the case $i=r, j\neq r$ and $j<r$.
Similarly one can show that the commutator relation holds for the case $i=r, j\neq r$ and $j>r$
and the case $i\neq r, j=r$. We now show the remaining case $i=j=r$.  In this case, a direct simplification yields 
\begin{align}
\label{commutator}
\begin{split}
e_i f_i- f_i e_i =
(\E_i \E_{i+1} \F_{i+1} \F_i - \F_{i+1} \F_i \E_i \E_{i+1} )
- v^{\ve} ( \E_{i+1} \F_{i+1} \E_i \F_i - \F_{i+1} \F_i \E_{i+1} \E_i) \\
+ v^{\ve} ( \E_i \F_i \E_{i+1} \F_{i+1} - \F_i \F_{i+1} \E_{i} \E_{i+1} ) 
+ ( \E_{i+1}\E_i \F_i \E_{i+1} - \F_i \F_{i+1} \E_{i+1} \E_i).
\end{split}
\end{align}
By using the commutator relations in $\U(\hsl_{n+1})$ on $\E_{i+1}\F_{i+1}$ and $\F_i\E_i$, we see that the term in the first parenthesis in (\ref{commutator}) is equal to
\begin{align}
\label{comm-1}
\E_i \frac{\K_{i+1} - \K^{-1}_{i+1}}{v-v^{-1}} \F_i + \F_{i+1} \frac{\K_i - \K^{-1}_i}{v-v^{-1}} \E_{i+1} . 
\end{align}
Note that $\F_{i+1} \F_i \E_{i+1} \E_i= \F_{i+1} \E_{\i+1} \F_i \E_i$, and by applying the commutator relation on 
both $ \F_{i+1} \E_{\i+1} $ and $\F_i \E_i$, 
the term in the second parenthesis in (\ref{commutator})  is simplified to
\begin{align}
\label{comm-2}
-v^{\ve} \left (
\frac{\K_{i+1} - \K^{-1}_{i+1}}{v-v^{-1}} \E_i \F_i 
+ \E_{i+1} \F_{i+1} \frac{\K_i-\K^{-1}_i}{v-v^{-1}} 
- \frac{\K_i-\K^{-1}_i}{v-v^{-1}} 
\frac{\K_{i+1} - \K^{-1}_{i+1}}{v-v^{-1}}
\right ).
\end{align}
By doing the same operation on the second monomial in 
the term in the third parenthesis in (\ref{commutator}), the term is simplied to
\begin{align}
\label{comm-3}
-v^{-\ve} \left (
\frac{\K_i -\K^{-1}_i}{v-v^{-1}} \E_{i+1} \F_{i+1} + 
\E_i \F_i \frac{\K_{i+1} - \K_{i+1}^{-1}}{v-v^{-1}} 
- \frac{\K_i- \K^{-1}_i}{v-v^{-1}} \frac{\K_{i+1}-\K^{-1}_{i+1}}{v-v^{-1}}
\right ).
\end{align}
Just like the term in the first parenthesis, the term in the last parenthesis in (\ref{commutator}) is simplified to 
\begin{align}
\label{comm-4}
\E_{i+1} \frac{\K_i-\K^{-1}_i}{v-v^{-1}} \F_{i+1} 
+\F_i \frac{\K_{i+1} -\K^{-1}_{i+1}}{v-v^{-1}} \E_i.
\end{align}
The sum of the terms in (\ref{comm-1})-(\ref{comm-4}) having $\E_i\F_i$ or $\F_i\E_i$ can be simplified to
\begin{align}
\label{comm-5}
\frac{-v^{-1} \K_{i+1} + v \K^{-1}_{i+1}}{v-v^{-1}} \cdot  \frac{\K_i-\K^{-1}_i}{v-v^{-1}}.
\end{align}
The sum of the terms in (\ref{comm-1})-(\ref{comm-4}) having $\E_{i+1}\F_{i+1}$ or $\F_{i+1}\E_{i+1}$ can be simplified to
\begin{align}
\label{comm-6}
\frac{-v^{-1} \K_i + v \K^{-1}_i}{v-v^{-1}} \cdot \frac{\K_{i+1} -\K^{-1}_{i+1}}{v-v^{-1}}.
\end{align}
The term involving only $\K_i^{\pm 1}$ and $\K^{\pm 1}_{i+1}$ in (\ref{comm-2}) plus the term in (\ref{comm-4}) is equal to
\begin{align}
\label{comm-7}
\K^{-1}_{i+1} \frac{\K_i -\K^{-1}_i}{v-v^{-1}}.
\end{align}
The term involving only $\K_i^{\pm 1}$ and $\K^{\pm 1}_{i+1}$ in (\ref{comm-3}) plus the term in (\ref{comm-5}) is equal to
\begin{align}
\label{comm-8}
\K_i \frac{\K_{i+1} - \K^{-1}_{i+1}}{v-v^{-1}}.
\end{align} 
So the commutator $e_if_i -f_i e_i$ for $i=r$ is equal to the sum of (\ref{comm-7}) and (\ref{comm-8}), which is
\[
\K^{-1}_{i+1} \frac{\K_i -\K^{-1}_i}{v-v^{-1}}+ \K_i \frac{\K_{i+1} - \K^{-1}_{i+1}}{v-v^{-1}}
=\frac{\K_i \K_{i+1} - \K^{-1}_i \K^{-1}_{i+1}}{v-v^{-1}} = \frac{k_i-k^{-1}_i}{v-v^{-1}}. 
\]
The commutator relation  (\ref{R3}) for the case $i=j=r$ is proved. 
This finishes the proof that the triple $(e_i, f_i, k^{\pm 1}_i)$ satisfies the relation  (\ref{R3}).

Next, we show that the pair $(e_i, e_j)$ satisfies the quantum Serre relation (\ref{R4}). 
First we observe that all cases are a consequence of the corresponding defining relations of $\U(\hsl_{n+1})$, except the cases
$(r\pm1, r)$ and  $(r, r\pm1)$.
For the case $(r-1, r)$, we can simplify the quantum Serre relation, say $S_{r-1, r}$,  as follows.
\begin{align*}
&\E^{(2)}_{r-1} (\E_r \E_{r+1} - v^{\ve} \E_{r+1} \E_r) -
\E_{r-1} (\E_r \E_{r+1} - v^{\ve} \E_{r+1} \E_r) \E_{r-1} 
+ (\E_r \E_{r+1} - v^{\ve} \E_{r+1} \E_r)  \E^{(2)}_{r-1} \\
 &=
\E^{(2)}_{r-1} \E_r\E_{r+1} - \E_{r-1} \E_{r} \E_{r+1} \E_{r-1} + \E_r \E_{r+1} \E^{(2)}_{r-1}\\
& \hspace{5cm} -v^{\ve}(
\E^{(2)}_{r-1} \E_{r+1} \E_r - \E_{r-1} \E_{r+1} \E_r\E_{r-1} + \E_{r+1} \E_r\E^{(2)}_{r-1}
)\\
& =
(\E^{(2)}_{r-1} \E_r -\E_{r-1} \E_r\E_{r-1} + \E_r \E^{(2)}_{r-1}) \E_{r+1}
- v^{\ve} \E_{r+1} (
\E^{(2)}_{r-1} \E_r - \E_{r-1} \E_r \E_{r-1} + \E_r \E^{(2)}_{r-1}
)\\
& =0.
\end{align*}
So the quantum Serre relation holds for the case $(i, j)=(r-1, r)$. The case for $(i, j)=(r+1, r)$ can be shown in exactly the same way.
Now we show the quantum Serre relation for the case $(r, r-1)$.
First of all, we need the formula (\ref{e-p}) for $p=2$, which reads
\[
e^{(2)}_r= \E^{(2)}_r \E^{(2)}_{r+1} - v^{\ve} \E_{r+1} \E^{(2)}_r \E_{r+1} + v^{2\ve} \E^{(2)}_{r+1} \E^{(2)}_{r}.
\]
By using this formula, we see immediately that the left-hand side of the quantum Serre relation 
for $(r, r-1)$, i.e.,
$e^{(2)}_r e_{r-1} - e_{r-1} e_r e_{r-1} + e_{r-1} e^{(2)}_r$,
is equal to 
\begin{align}
\label{Serre}
\begin{split}
&
(\E^{(2)}_r \E^{(2)}_{r+1} \E_{r-1} - \E_r \E_{r+1} \E_{r-1} \E_r\E_{r+1} + \E_{r-1} \E^{(2)}_r \E^{(2)}_{r+1} )\\
&-v^{\ve}(
\E_{r+1} \E^{(2)}_r \E_{r+1} \E_{r-1} 
-\E_{r+1} \E_r \E_{r-1} \E_{r+1} - \E_r \E_{r+1} \E_{r-1} \E_{r+1} \E_r
+ \E_{r-1} \E_{r+1} \E^{(2)}_r \E_{r+1}
)\\
&+v^{2\ve} (
E^{(2)}_{r+1} \E^{(2)}_r \E_{r-1} 
-\E_{r+1} \E_r \E_{r-1} \E_{r+1}\E_r
+ \E_{r-1} \E^{(2)}_{r+1}\E^{(2)}_r
).
\end{split}
\end{align}
Replace the piece $\E_{r+1} \E_r \E_{r+1}$ by $\E_{r} \E_{r+1}^{(2)} + \E^{(2)}_{r+1} \E_r$ in the second monomial, we see that the terms
in the first parenthesis in (\ref{Serre}) are equal to 
\begin{align}
\label{Serre-1}
- \E_r\E_{r-1} \E^{(2)}_{r+1} \E_r.
\end{align}
Similarly, the third term with the coefficient $v^{2\ve}$ in (\ref{Serre}) is equal to 
\begin{align}
\label{Serre-2}
- v^{2\ve} \E_r \E^{(2)}_{r+1} \E_{r-1} \E_r.
\end{align}
The monomials in (\ref{Serre-1}) and (\ref{Serre-2}) add up to 
\begin{align}
\label{Serre-3}
-v^{\ve} [2] \E_r \E_{r-1} \E_{r+1}^{(2)} \E_r.
\end{align}
The first, second and fourth monomials in the term with coefficient $v^{\ve}$ add up to zero. 
The remaining monomial  in the term with $v^{\ve}$ is
$v^{\ve} \E_r \E_{r+1} \E_{r-1} \E_{r+1} \E_r = v^{\ve} [2] \E_r \E_{r-1} \E_{r+1}^{(2)} \E_r$, which cancels with (\ref{Serre-3}).
So the quantum Serre relation for $(r, r-1)$ holds. 
This finishes the proof of the quantum Serre relation (\ref{R4}) for the pair $(e_i, e_j)$.

Similarly one can show that the quantum Serre relation (\ref{R5}) for the pair $(f_i, f_j)$ holds.
One can also apply the involution  on $\U(\hsl_{n+1})$, 
defined by $\E_i\mapsto \F_i, \F_i\mapsto \E_i, \K_i\mapsto \K^{-1}_i$ for all $i\in I_{n+1}$, to obtain the proof.

The proof of $\Phi_{r,\ve}$ being an algebra homomorphism for $n>2$ is thus finished. 

\subsection{Existence of $\Phi_{r,\ve}$ in the $n=2$ case}
\label{Proof-2}
In this section, we  assume that $n=2$, and we shall show that the set $\{ e_i, f_i, k_i^{\pm 1} | i\in I_2\}$ satisfies the defining relations
of $\U(\hsl_2)$. 
The relations (\ref{R1}) and (\ref{R2}) can be verified directly, just like the $n>2$ case, with very minor modifications. 
The commutator relation (\ref{R3}) can be proved in exactly the same way as the $n>2$ case.
So we only need to verify the quantum Serre relations (\ref{R4}) and (\ref{R5}). 

Let us verify the relation (\ref{R4}).
For simplicity, we use $(i,j, \ell)$ for $(r, r+1, r-1)$. 
The first relation to show is 
\begin{align}
\label{2-Serre}
e^{(3)}_\ell e_i- e^{(2)}_\ell e_i e_\ell + e_\ell e_i e^{(2)}_\ell - e_i e^{(3)}_\ell=0.
\end{align}
A simplification yields that the left-hand side of (\ref{2-Serre}) is equal to
\begin{align}
\label{2-Serre-1}
\begin{split}
(\E^{(3)}_\ell \E_i \E_j - \E^{(2)}_\ell \E_i\E_j \E_\ell + \E_\ell \E_i \E_j \E^{(2)}_\ell - \E_i\E_j \E^{(3)}_\ell)\\
-v^{\ve}(
\E^{(3)}_\ell \E_j \E_i -\E^{(2)}_\ell \E_j\E_i \E_\ell + \E_l \E_j\E_i \E^{(2)}_\ell - \E_j \E_i \E^{(2)}_\ell
).
\end{split}
\end{align}
Observe that if the term in the first parenthesis is zero, so is the term in the second one by switching the roles of $i$ and $j$.
Recall that we have
\[
\E^{(3)}_\ell \E_i  = \E_\ell \E_i \E^{(2)}_\ell -[2] \E_i \E^{(3)}_\ell,
\E_i \E^{(3)}_\ell = \E^{(2)}_\ell \E_i \E_\ell - [2] \E^{(3)}_\ell \E_i.
\]
Applying the above identities to the first and fourth monomials in the term in the first parenthesis   in (\ref{2-Serre-1}), we have
\begin{align}
\label{Serre-deg-5}
\begin{split}
&\E^{(3)}_\ell \E_i \E_j - \E^{(2)}_\ell \E_i\E_j \E_\ell + \E_\ell \E_i \E_j \E^{(2)}_\ell - \E_i\E_j \E^{(3)}_\ell\\
&=\E_\ell \E_i (\E^{(2)}_\ell \E_j + E_j\E^{(2)}_\ell) -
(\E^{(2)}_\ell \E_i + \E_i \E^{(2)}_\ell) \E_j \E_\ell\\
&= \E_\ell \E_i \E_\ell \E_j\E_\ell -  \E_\ell \E_i \E_\ell \E_j\E_\ell=0.
\end{split}
\end{align}
Therefore, the relation (\ref{2-Serre}) holds.

The second relation we need to verify is 
\begin{align}
\label{2-Serre-b}
S_{i\ell}\equiv e^{(3)}_i e_\ell - e^{(2)}_i e_\ell e_i + e_i e_\ell e^{(2)}_i - e_\ell e^{(3)}_i =0.
\end{align}
Instead of proving it directly, which involves monomials of degree 7, 
we shall make use of Lusztig's twisted derivation $r_k$ on the positive half of a quantum group in~\cite[1.2.13]{L10}. 
Let $\U^+(\hsl_3)$ be the positive half of $\U(\hsl_3)$ generated by $\E_k$ for all $k\in I_3$. 
To each $k$, there is a linear map $r_k: \U^+(\hsl_3) \to \U^+(\hsl_3)$ such that 
\[
r_k (1) =0, r_k (\E_{k'})=\delta_{k, k'}, r_k(xy) = v^{k \cdot |y|} r_k(x) y + x r_k(y),
\]
for all homogeneous element $y$. 
Here $|y|=(y_0, y_1, y_2)\in \mbb Z^{I_3}$ is the degree of $y$ and $k\cdot |y| = 2y_k - y_{k+1} - y_{k-1}$. 
Recall $(r, r+1, r-1)=(i,j,\ell)$. 
For now we assume that $\ve =-1$. 
We will freely use the following formulas in later analysis. 

\begin{lem}
For all $p\geq 0$, we have 
\begin{align}
%
\label{s1}
\tag{s1}
r_\ell (e_i^{(p)}) &=0,\\
\label{s2}
\tag{s2}
r_i(e^{(p)}_i)& =0, \\
\label{s3}
\tag{s3}
r_j(e^{(p)}_i) &= (v-v^{-1}) v^{p-2} \sum_{a=0}^{p-1}  (-1)^a v^{-2a} \E^{(a)}_j \E^{(p)}_i \E^{(p-1-a)}_j,\\
\label{s4}
\tag{s4}
r_i r_j (e^{(p)}_i) & = (v-v^{-1})v^{p-2} e^{(p-1)}_i,\\
\label{s5}
\tag{s5}
r_jr_j (e^{(p)}_i) & = (v-v^{-1}) (v^2-v^{-2}) v^{2(p-3)} \sum_{a=0}^{p-2} v^{-3a} \E^{(a)}_j \E^{(p)}_i \E^{(p-2-a)}_j.
\end{align}
\end{lem}

\begin{proof}
The equality (\ref{s1}) is because the $\ell$-degree of $e_i$ is zero.
For (\ref{s2}), we only need to show that $r_i(e_i) =0$, which can be done as follows.
\[
r_i(e_i) = r_i(\E_i\E_j - v^{-1} \E_j\E_i) = v^{i\cdot j} r_i(\E_i) \E_j - v^{-1} \E_j r_i(\E_i) = v^{-1} \E_j - v^{-1} \E_j=0. 
\]
We show (\ref{s3}) by induction. When $p=1$, we have
\[
r_j(e_i) =  r_j(\E_i\E_j - v^{-1} \E_j\E_i)= \E_i - v^{-1} v^{i\cdot j} r_j(\E_j) \E_i  = (v-v^{-1}) v^{-1} \E_i. 
\]
Assume that (\ref{s3}) holds for $p$, we want to show that it holds for $p+1$.
We have 
\begin{align}
\label{rji}
r_j(e^{(p+1)}_i) =\frac{1}{[p+1]} r_j (e^{(p)}_i . e_i) 
=\frac{1}{[p+1]}\left (v r_j(e^{(p)}_i) e_i + e^{(p)}_i r_j(e_i) \right )
\end{align}
We simplify the first term as follows.
\begin{align}
\label{rji-1}
\begin{split}
&\frac{1}{[p+1]} v r_j(e^{(p)}_i) e_i
= \frac{(v-v^{-1})v^{p-1} }{[p+1]} \sum_{a=0}^{p-1} (-1)^a v^{-2a} \E^{(a)}_j \E^{(p)}_i \E^{(p-1-a)}_j (\E_i\E_j-v^{-1}\E_j\E_i) \\
& =(v-v^{-1}) v^{p-1} \sum_{a=0}^{p-1} (-1)^a v^{-2a} \E^{(a)}_j \E_i^{(p+1)} \E_j^{(p-a)} \\
& +
\frac{(v-v^{-1}) v^{p-1}}{[p+1]} \sum_{a=0}^{p-1} (-1)^a (v^{-2a} [p-1-a] - v^{-1-2a}[p-a]) \E^{(a)}_j \E^{(p)}_i \E^{(p-1}_j \E_i,
\end{split}
\end{align}
where we use the fact 
\[
\E^{(p-1-a)}_j \E_i \E_j = \E_i \E^{(p-a)}_j + [p-1-a] \E^{(p-a)}_j \E_i.
\]
The second term can be simplified as follows.
\begin{align}
\label{rji-2}
\begin{split}
\frac{1}{[p+1]} e^{(p)}_i r_j (e_i) =
(v-v^{-1}) v^{p-1} (-1)^{p} v^{-2p} \E^{(p)}_j \E^{(p+1)}_i \\ +
\frac{(v-v^{-1})v^{-1}} {[p+1]} \sum_{a=0}^{p-1} (-1)^a v^{-a} \E^{(a)}_j \E^{(p)}_i \E^{(p-a)}_i \E_i. 
\end{split}
\end{align}
The second terms  in (\ref{rji-1}) and (\ref{rji-2}) cancel once we observe that
\[
v^{p-2a}[p-1-a] - v^{p-1-2a}[p-a] + v^{-a}=0.
\] 
The equality (\ref{s3}) is then followed by adding the first terms in (\ref{rji-1}) and (\ref{rji-2}).
Now we show (\ref{s4}). By (\ref{s3}), we have
\begin{align}
\begin{split}
r_i r_j (e^{(p)}_i) &=(v-v^{-1}) v^{p-2} \sum_{a=0}^{p-1}  (-1)^a v^{-2a} \E^{(a)}_j   v^{-(p-1- a)} r_i(\E^{(p)}_i) \E^{(p-1-a)}_j\\
&=(v-v^{-1}) v^{p-2} \sum_{a=0}^{p-1}  (-1)^a v^{-a} \E^{(a)}_j    \E^{(p-1)}_i \E^{(p-1-a)}_j\\
&=(v-v^{-1}) v^{p-2} e^{(p-1)}_i.
\end{split}
\end{align}
So (\ref{s4}) holds. 
Finally we show (\ref{s4}). By (\ref{s3}), we have
\begin{align}
\begin{split}
& r_j r_j (e^{(p)}_i) \\
&= (v-v^{-1}) v^{p-2} \sum_{a=0}^{p-1}  (-1)^a v^{-2a} 
( v^{ p-2a-2} r_j(\E^{(a)}_j ) \E^{(p)}_i \E^{(p-1-a)}_j+ \E^{(a)}_j) \E^{(p)}_i r_j( \E^{(p-1-a)}_j))\\
& =(v-v^{-1}) v^{p-2} \sum_{a=0}^{p-1}  (-1)^a v^{-2a} 
( v^{ p-a-3} \E^{(a-1)}_j  \E^{(p)}_i \E^{(p-1-a)}_j+ v^{p-2-a} \E^{(a)}_j) \E^{(p)}_i \E^{(p-2-a)}_j)\\
&= (v-v^{-1}) v^{p-2} \sum_{a=0}^{p-2} ( -1)^a ( v^{p-3a-6} - v^{p-3a-2}) \E^{(a)}_j \E^{(p)}_i \E^{(p-a-2)}_j\\
&=(v-v^{-1}) (v^2-v^{-2}) v^{2(p-3)}  \sum_{a=0}^{p-2} ( -1)^a v^{-3a}  \E^{(a)}_j \E^{(p)}_i \E^{(p-a-2)}_j.
\end{split}
\end{align}
So (\ref{s5}) holds. 
The lemma is therefore proved.
\end{proof}

%
%
%
%
%
%
%
%
%
%
%
  
In light of~\cite[Lemma 1.2.15]{L10},  to show that $S_{i\ell}=0$, it is sufficient to show that 
\begin{align}
\label{Sil-condition}
r_i(S_{i\ell}) =0, r_j(S_{i\ell}) =0, r_\ell (S_{i,\ell})=0. 
\end{align}
Due to (\ref{s2}) and $r_i(\E_\ell)=0$, we have immediately that $r_i(S_{i\ell})=0$. 
Moreover, we have
\begin{align*}
r_\ell (S_{i\ell})
&=e^{(3)}_i - e^{(2)}_i r_\ell (e_\ell e_i) + e_i r_\ell ( e_\ell e^{(2)}_i) - 
v^{\ell \cdot 3(i+j)} e^{(3)}_i\\
& = e^{(3)}_i - v^{-2} e^{(2)} e_i + v^{-4} e_i e^{(2)}_i -v^{-6} e^{(3)}_i\\
&=(1 - v^{-2}[3]+ v^{-4} [3] - v^{-6}) e^{(3)}_i =0. 
\end{align*}
So it remains to show that  $r_j(S_{i\ell})=0$.
This is further reduced to show that 
\begin{align}
\label{rj-condition}
r_i r_j(S_{i\ell})=0, r_j r_j(S_{i\ell})=0, r_\ell r_j(S_{i\ell})=0.
\end{align}
By a direct computation, we get
\begin{align}
\label{rjS}
\begin{split}
r_j(S_{i\ell}) =
v^{-1} r_j(e^{(3)}_i) \E_\ell - r_j (e^{(2)}_i) \E_\ell e_i - v^{-1}(v-v^{-1}) e^{(2)}_i \E_\ell \E_i \\
+ (v-v^{-1}) \E_i\E_\ell e^{(2)}_i + e_i \E_\ell r_j (e^{(2)}_i) -\E_\ell r_j (e^{(3)}_i). 
\end{split}
\end{align}
By applying $r_i$ to the formula (\ref{rjS}) and  a direction computation, we have
\begin{align*}
&r_i r_j(S_{i\ell})  = v^{-2} r_i r_j (e^{(3)}_i) \E_\ell 
- (v-v^{-1}) e_i \E_\ell e_i \\
& - v^{-1}(v-v^{-1} )e^{(2)}_i \E_\ell
+ v(v-v^{-1}) \E_\ell e^{(2)}_i 
+ (v-v^{-1}) e_i \E_\ell e_i 
- \E_\ell r_ir_j(e^{(3)}_i) \\
&=
v^{-1}(v-v^{-1}) e^{(2)}_i \E_\ell  
- v^{-1}(v-v^{-1} ) e^{(2)}_i \E_\ell + 
v(v-v^{-1}) \E_\ell e^{(2)}_i - v(v-v^{-1}) \E_\ell e^{(2)}_i 
=0 .
\end{align*}
So to show that $r_j(S_{i\ell})=0$, it remains to show that $r_\ell r_j (S_{i\ell}) =0$ and $r_jr_j(S_{i\ell})=0$.
By a direct computation, we get
\begin{align*}
 r_\ell r_j (S_{i\ell}) =
v^{-3} (v^{2}-v^{-2}) r_j (e^{(3)}_i)  - v^{-2} r_j (e^{(2)}_i) e_i - v^{-2}(v-v^{-1}) e^{(2)}_i \E_i \\
+ v^{-4}(v-v^{-1}) \E_i e^{(2)}_i + v^{-3} e_i r_j(e^{(2)}_i) .
\end{align*}
Since there is no $\E_\ell$ in $r_\ell r_j (S_{i\ell}) $, so $r_\ell r_\ell r_j(S_{i\ell})=0$. 
Next, we compute $r_i r_\ell r_j (S_{i\ell})$ as follows.
\begin{align*}
r_i r_\ell r_j (S_{i\ell})
=
[v^{-2} (v^2-v^{-2}) (v-v^{-1}) 
- v^{-1} (v-v^{-1}) [2]  
+ v^{-3} (v-v^{-1}) [2]] e^{(2)}_i
=0. 
\end{align*}
Further we compute $r_jr_\ell r_j(S_{i\ell})$. We have 
\begin{align*}
r_j r_\ell r_j (S_{i\ell}) =
v^{-3}(v-v^{-1})^3
(\E^{(3)}_i \E_j - v^{-1} \E^{(2)}_i \E_j \E_i + v^{-2} \E_i \E_j \E^{(2)}_i - v^{-3} \E_j \E^{(3)}_i)=0.
\end{align*}
The above analysis shows that $ r_\ell r_j (S_{i\ell})=0$.  
So to show that $r_j (S_{i \ell})=0$ it remains to show that $r_j r_j (S_{i\ell})=0$. 
By a direct computation, we have
\begin{align}
\label{rjj}
\begin{split}
\frac{r^2_j(S_{i\ell}) } { v^{-2}(v-v^{-1}) (v^2- v^{-2}) } &
= (\E^{(3)}_i \E_j - v^{-3} \E_j \E^{(3)}_i) \E_\ell 
-\E^{(2)}_i \E_\ell e_i \\
&\hspace{-1cm} -v^{-1} (v \E^{(2)}_i \E_j - v^{-1}\E_j \E^{(2)}_i) \E_\ell \E_i
+ \E_i \E_\ell (v \E^{(2)}_i \E_j - v^{-1} \E_j \E^{(2)}_i) \\
&\hspace{3.5cm} + e_i \E_\ell \E^{(2)}_i - v^2 \E_\ell( \E^{(3)}_i \E_j - v^{-3} \E_j \E_i^{(3)}) .
\end{split}
\end{align}
Now substitute $e_i$ by $\E_i\E_{i+1}-v^{-1}\E_{i+1}\E_i$, we see that 
(\ref{rjj}) is equal to 
\begin{align}
\label{rjj-1}
\begin{split}
\E^{(3)}_i \E_j\E_\ell - \E^{(2)}_i \E_j \E_\ell \E_i + \E_i \E_j \E_\ell \E^{(2)}_i \\
-v^{-3} \E_j \E^{(3)} \E_\ell + v^{-2} \E_\ell \E^{(2)}_i \E_\ell \E_i - v^{-1} \E_j\E_i \E_\ell \E^{(2)}_i \\
-\E^{(2)}_i \E_\ell \E_i \E_j + v \E_i \E_\ell \E^{(2)}_i \E_j - v^{2} \E_\ell \E^{(3)}_i \E_j \\
+\E^{(2)}_i \E_\ell \E_j \E_i - v^{-1} \E_i \E_\ell \E_j \E^{(2)}_i + v^{-1} \E_\ell \E_j \E^{(3)}_i.
\end{split}
\end{align}
By (\ref{Serre-deg-5}), the first row in (\ref{rjj-1}) is equal to $\E_j \E_{\ell} \E^{(3)}_i$, which cancels with the second row.
The third row in (\ref{rjj-1}) is equal to $-v^{-1} \E^{(3)}_i \E_{\ell} \E_j$, which cancels with the fourth row.
So we get $r^2_j(S_{i\ell})=0$. 
This finishes the proof of (\ref{rj-condition}), and therefore $r_j(S_{i\ell})=0$. 
In turn, this shows that (\ref{Sil-condition}) holds, and thus $S_{i\ell}=0$, i.e., (\ref{2-Serre-b}), as desired. 
The above proof assumes that $\ve=-1$. The case for $\ve=1$ can be proved by rewriting $e_i$ 
as $v^\ve (\E_j \E_i - v^{-\ve} \E_i\E_j)$ and the proof for $\ve=1$ case applies by switching the roles of $i$ and $j$. 

The relation (\ref{R5}) is a consequence of the relation (\ref{R4}) by applying the involution on $\U(\hsl_{3})$ 
sending $\E_k$ to $\F_k$ $\forall k\in I_3$. 
This finishes the proof of  $\Phi_{r,\ve}$ being an algebra homomorphism for $n=2$.

\subsection{Proof of Proposition~\ref{variant}}
\label{Proof-variant}
In this section, we provide a proof of Proposition~\ref{variant}.
The last two relations in (\ref{gl-def}) have been verified in the proof of Theorem~\ref{Phi}. 
The relations in the first three rows of (\ref{gl-def}) can be checked directly. 
By using the relation $\prod_{i=1}^{n+1} \breve K_i =1$, we see that
$k_i = l_i l^{-1}_{i+1}$ for all $i\in I_n$. Indeed, the equality holds for all $i\in I_{n}-\{ n\}$ obviously.
If $i=n$, then 
\begin{align*}
l_n l^{-1}_{n+1} =l_n l^{-1}_1 = \breve K^{-1}_{r+1} \cdots \breve K^{-1}_n  (\breve K_1\cdots \breve K_r)^{-1} = \breve K_{n+1} = k_n. 
\end{align*} 
So the commutator relation in (\ref{gl-def}) is a consequence of the commutator relation of $\U(\hsl_{n+1})$. 
This implies that the map $\Phi_{r,\ve: 1}$ is an algebra homomorphism. The injective property follows from the triangular decompositions of $\U(\hgl_n)$ and $\U(\hsl_{n+1})_1$. 
This finishes the proof.

\end{document}